\documentclass[preprint,12pt]{elsarticle}

\usepackage{amssymb}
\usepackage{amsmath}
\usepackage{mathtools}
\usepackage{amsthm}
\usepackage{color}
\usepackage{multirow}
\usepackage{lscape}

\usepackage{subfigure}
\usepackage{float}

\usepackage[svgnames]{xcolor}
\usepackage{soul}

\usepackage{algorithm2e}  

\usepackage{tikz}
\tikzset{every label/.style={font=\footnotesize,inner sep=1pt}}

\usetikzlibrary{graphs,graphs.standard,shapes,arrows.meta,chains}

\usepackage[font=small,labelfont=bf]{caption}

\usepackage{hyperref}

\newcommand{\kmin}{\kappa_{\rm min}}
\newcommand{\kmax}{\kappa_{\rm max}}
\journal{Journal of Computational Physics}

\begin{document}

\begin{frontmatter}



\title{Speeding up a few orders of magnitude the Jacobi method: high order Chebyshev-Jacobi over GPUs}


\author{J.E.~Adsuara\corref{cor1}\fnref{label1}} 
\author{M.A.~Aloy\corref{cor1}\fnref{label1}}
\author{P.~Cerd\'a-Dur\'an\corref{cor1}\fnref{label1}}
\author{I.~Cordero-Carri\'on\corref{cor1}\fnref{label2}}
\cortext[cor1]{jose.adsuara@uv.es, miguel.a.aloy@uv.es, pablo.cerda@uv.es, isabel.cordero@uv.es}
\address[label1]{Departamento de Astronom\'{\i}a y Astrof\'{\i}sica, Universidad de Valencia, E-46100, Burjassot, Spain.}
\address[label2]{Departamento de Matem\'aticas, Universidad de Valencia, E-46100, Burjassot, Spain.}

\begin{abstract}
In this technical note we show how to reach a remarkable speed up when solving elliptic partial differential equations with finite differences thanks to the joint use of the Chebyshev-Jacobi method with high order discretizations and its parallel implementation over GPUs.
\end{abstract}


\end{frontmatter}


\section{Introduction}
\label{sec:int}
In \cite{YangMittal2014, Adsuetal15, Adsuara2016369} a generalization
of the weighted Jacobi method for linear systems was
presented and improved. Adsuara et al.
\cite{Adsuara2017446} showed that this new Chebyshev-Jacobi
  method (CJM) is equivalent to a Generalized Richardson iterative
  scheme. These authors showed how to obtain analytically the weights
not only for the classical 5-points discretization of the Laplacian
operator, but also for high order discretizations, involving a larger number of points. They found that using higher-order discretizations of the
  Laplacian is always advantageous, since they lead to a solution of
the same quality than lower-order discretizations, but in
a shorter time. This point in combination with the trivial
parallelization of the method predicts a promising
performance. Here, we go an step forward by parallelizing the CJM for
GPUs.

This note is organized as follows. In Section \ref{sec:met} we
recap the essentials
  of the CJM with high-order discretizations and the
  methodological steps to port it to the GPU technology. In Section
\ref{sec:res} we show
  the speed up of the new implementation in various GPU
  architectures. Finally, in the last section, we summarize the main
achievements of the article and define the outlines for the future
work.

\section {Methodology}
\label{sec:met}

In the following subsections we first present a brief
  summary of the CJM and then comment on the various ingredients involved in our analysis of the
acceleration factors for the different GPU architectures.

\subsection{The CJM with a high order discretization of the Laplacian}
\label{subsec:cjHigOrd}

 As showed in \mbox{\cite{Adsuara2017446}},
the CJM is the optimal one of all the possible
schemes introduced previously in \cite{YangMittal2014, Adsuetal15,
  Adsuara2016369}. Briefly, the CJM consists on a weighted classical
Jacobi method with strictly different weights at each iteration. A
particular scheme is defined then by a set of weights, which is
precalculated before starting the iterative algorithm as a
transformation of the zeros of a Chebyshev polynomial. This
transformation depends, among other aspects, on the resolution of the
mesh, the boundary conditions and the required
tolerance. The grid resolution and boundary conditions
  determine the wave numbers of elementary trigonometric
  waves that one can fit on the mesh. This allow us to bound the eigenvalues 
  of the iteration matrix to the interval $[\kmin,\kmax]$, whose boundaries can
  be obtained analytically (see below).

When finite differences are used in order to solve elliptic partial differential equations, the problem reduces to solving a linear system. A CJM can be used and the corresponding transformation of the zeros of a Chebyshev polynomial also depends on the resulting matrix associated to the linear system. As it has been already mentioned in \cite{Adsuara2017446}, the classical 5-points discretization of the bidimensional Laplacian leads to a consistently ordered (CO) matrix, where Young's theory is applicable. However, there are many other useful discretizations \cite{Adams:1988}, which end up with a non-consistently ordered matrix, where Young's theory cannot be used.

In \cite{Adsuara2017446} we presented two
parametrized discretizations of the
  Laplacian in two spatial dimensions, which can achieve an order of accuracy larger than two. In this note
we use one discretization of each of these two families. From the
first family, we use  a 9-point discretization corresponding to
a value of the parameter in \cite{Adsuara2017446} of $ \alpha=\frac{2}{3}$, whose
Laplacian is discretized as
\begin{gather}
\Delta u_{ij}=\frac{1}{6 h^2} \bigg [ 4 u_{i-1,j} + 4 u_{i+1,j} + 4u_{i,j-1} + 4 u_{i,j+1} \nonumber \\
+ u_{i-1,j-1} + u_{i+1,j+1} + u_{i-1,j+1} + u_{i+1,j-1} - 20 u_{i,j} \bigg ]\, 
\label{eq:9-points}
\end{gather}  
and whose $\kmin$ and $\kmax$ values, needed for obtaining the scheme and related with the size of the mesh, are
\begin{gather}
\kmin^{(9)} = \frac{4}{5} \bigg [ \sin^2{\frac{\pi}{2 N_x}} + \sin^2{\frac{\pi}{2 N_y}} \bigg ] \\
+ \frac{1}{5} \bigg [ \sin^2{ \bigg ( \frac{\pi}{2 N_x} +\frac{\pi}{2 N_y} \bigg )} + \sin^2{ \bigg ( \frac{\pi}{2 N_x} - \frac{\pi}{2 N_y} \bigg )} \bigg ] ,\\
\kmax^{(9)} = \frac{8}{5} \,,
\label{Eq:kmkM9p}
\end{gather}  
where the subscript $(9)$ refers to the number of points
  in which the Laplacian operator is discretized, and $N_x$ and $N_y$
  are the number of points per dimension in the spatial directions $x$
  and $y$, respectively. The grid assumed in this discretization is
  uniform, so that $L_x/N_x=L_y/N_y:=h$, $L_x$ and $L_y$ being the
  sizes along the $x-$ and $y-$ directions, respectively.

From the second family of discretizations, we use  a
17-point discretization with $\alpha = \frac{2}{3}$, whose
 discrete representation of the Laplacian
 reads
\begin{gather}
\Delta u_{ij}  = \frac{1}{72 h^2} \bigg [
- 4 u_{i-2,j} + 64 u_{i-1,j} + 64 u_{i+1,j} - 4 u_{i+2,j} \nonumber\\
- 4 u_{i,j-2} + 64 u_{i,j-1} + 64 u_{i,j+1} - 4 u_{i,j+2}  \nonumber\\
- u_{i-2,j-2} + 16 u_{i-1,j-1} + 16  u_{i+1,j+1} - u_{i+2,j+2} \nonumber\\ 
- u_{i-2,j+2} + 16 u_{i-1,j+1} + 16  u_{i+1,j-1} - u_{i+2,j-2} - 300 u_{i,j} \bigg ].
\label{eq:17-points}
\end{gather}  
and whose  corresponding $\kmin$ and $\kmax$ values are
\begin{gather}
\kmin^{(17)} =
-\frac{4}{75} \bigg [ \sin^2 \frac{\pi }{N_x} + \sin^2 \frac{\pi }{N_y} \bigg ] + \frac{64}{75} \bigg [ \sin^2 \frac{\pi}{2
  N_x}  + \sin^2 \frac{\pi}{2 N_y} \bigg ] \nonumber\\
- \frac{1}{75} \bigg [ \sin^2 { \bigg ( \frac{\pi}{N_x} + \frac{\pi}{N_y} \bigg ) + \sin^2 { \bigg ( \frac{\pi}{N_x} - \frac{\pi}{N_y} \bigg ) }} \bigg ] \\
+ \frac{16}{75} \bigg [ \sin^2 { \bigg ( \frac{\pi}{2N_x} + \frac{\pi}{2N_y} \bigg ) + \sin^2 { \bigg ( \frac{\pi}{2N_x} - \frac{\pi}{2N_y} \bigg ) }} \bigg ] \\
\kmax^{(17)} = \frac{128}{75}.
\end{gather}  

\begin{figure*}
\centering
\begin{tikzpicture}[baseline={(current bounding box.center)},every path/.style={>=latex},every node/.style={draw,circle}]
  \node (i-1)         at (-1.5,0)  {$\frac{2}{3}$};
  \node (i-1--j+1)  at (-1.5,1.5)  {$\frac{1}{6}$};
  \node (i-1--j-1)   at (-1.5,-1.5) {$\frac{1}{6}$};
  \node (i+1--j-1)  at (1.5,-1.5)  {$\frac{1}{6}$};
  \node (i+1--j+1) at (1.5,1.5)    {$\frac{1}{6}$};
  \node (i)            at (0,0)    {$\frac{-10}{3}$};
  \node (i+1)        at (1.5,0)    {$\frac{2}{3}$};
  \node (j-1)         at (0,-1.5)   {$\frac{2}{3}$};
  \node (j+1)        at (0,1.5)    {$\frac{2}{3}$};
  \draw[very thick] 
        (j-1) -- (i)
        (i)   -- (j+1)
        (i-1) -- (i)
        (i)   -- (i+1)
        (i-1--j+1) -- (i)
        (i-1--j-1) -- (i)
        (i+1--j+1) -- (i)
        (i+1--j-1) -- (i);
\end{tikzpicture}
\qquad \qquad \qquad
\begin{tikzpicture}[baseline={(current bounding box.center)},every path/.style={>=latex},every node/.style={draw,circle}]
  \node (i-1)         at (-1.5,0)  {$\frac{8}{9}$};
  \node (i-1--j+1)  at (-1.5,1.5)  {$\frac{2}{9}$};
  \node (i-1--j-1)   at (-1.5,-1.5) {$\frac{2}{9}$};
  \node (i+1--j-1)  at (1.5,-1.5)  {$\frac{2}{9}$};
  \node (i+1--j+1) at (1.5,1.5)    {$\frac{2}{9}$};
  \node (i)            at (0,0)    {$\frac{-25}{6}$};
  \node (i+1)        at (1.5,0)    {$\frac{8}{9}$};
  \node (j-1)         at (0,-1.5)   {$\frac{8}{9}$};
  \node (j+1)        at (0,1.5)    {$\frac{8}{9}$};
  \node (i-2)         at (-3,0)   {$\frac{-1}{18}$};
  \node (i+2)        at (3,0)    {$\frac{-1}{18}$};
  \node (j-2)         at (0,-3)   {$\frac{-1}{18}$};
  \node (j+2)        at (0,3)    {$\frac{-1}{18}$};
  \node (i-2--j+2)  at (-3,3)   {$\frac{-1}{72}$};
  \node (i-2--j-2)   at (-3,-3)  {$\frac{-1}{72}$};
  \node (i+2--j-2)  at (3,-3)   {$\frac{-1}{72}$};
  \node (i+2--j+2) at (3,3)    {$\frac{-1}{72}$};
  \draw[very thick] 
        (j-1) -- (i)
        (i)   -- (j+1)
        (i-1) -- (i)
        (i)   -- (i+1)
        (j-2) -- (j-1)
        (j+1)   -- (j+2)
        (i-2) -- (i-1)
        (i+1)   -- (i+2)
        (i-1--j+1) -- (i)
        (i-1--j-1) -- (i)
        (i+1--j+1) -- (i)
        (i+1--j-1) -- (i)
        (i-2--j+2) -- (i-1--j+1)
        (i-2--j-2) -- (i-1--j-1)
        (i+2--j+2) -- (i+1--j+1)
        (i+2--j-2) -- (i+1--j-1);
\end{tikzpicture}
\caption{Schematic representation of the 9-points and 17-point
   discretizations of the Laplacian in two
    spatial dimensions, corresponding to a value
  $\alpha = \frac{2}{3}$. Each circle corresponds to a neighboring point of a given
    one $(i,j)$, represented by the central circle. The value enclosed
    by each circumference corresponds to the coefficient of the
    discretizations shown in Eqs.\,\ref{eq:9-points} and
    \ref{eq:17-points} multiplied by $h^2$. See \cite{Adsuara2017446}
  for details.}
\label{fig::fig01}
\end{figure*}
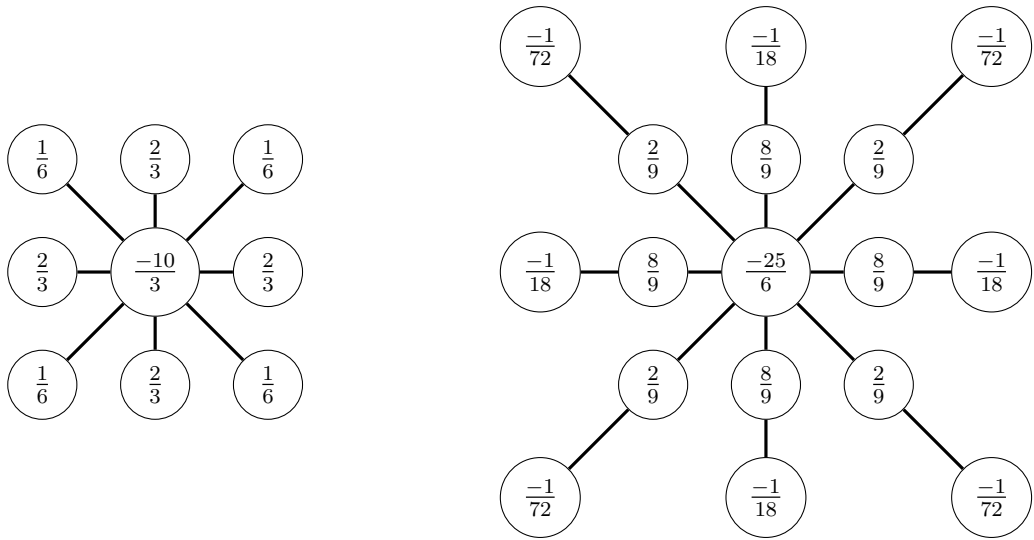

In Figure \ref{fig::fig01} we show an schematic view of the two
stencils for the two different discretizations along with the
numerical coefficients corresponding to each node.  Note that the
value of the discrete Laplacian operator evaluated at the central node
is a weighted sum over all the nodes schematically represented in each
of the figure panels. The purpose of these schemes is showing that
given a central point, the Laplacian discretization is fully specified
by providing a list of $M$ numerical coefficients enclosed in all
the surrounding nodes.

\subsection{GPUs}
\label{subsec:cheJacMet}

Modern
  GPUs have become extremely powerful processing units, which can
  accelerate enormously the computation of heavy mathematical
  operations. However, for that to happen, the algorithms must be
  properly vectorized and the pipeline of arithmetic operations and
  data transfers from/to the main memory must be optimized. The CJM method
  is perfectly suited for its porting to GPUs.  In each iteration of the
method the approximate solution is stored in each node
and we must apply the same set of elementary arithmetic operations to
each of them.

\section{Results}
\label{sec:res}

In this section we show the execution times obtained when solving a
particular problem in the two discretizations introduced above using
two specific GPUs models.

\subsection{CUDA devices}
\label{subsec:gpuDev}

We present in this subsection the specific devices we used. We aim to test more than a single GPU model in order
  to properly assess the scaling properties of our implementation of
  the CJM, which is based on the CUDA \cite{CUDA} technology of
NVDIA. CUDA devices
can be characterized by two parameters (among others):
the Compute Capability (CC) and the number of Streaming
Multiprocessors (SMs). The CC of a device is
represented by a version number. It comprises a major revision number
X and a minor revision number Y, and it is denoted with the format
``X.Y''. Devices with the same major revision number belong to the same core architecture.

We have used two CUDA devices: a Tesla K40c and a GeForce GTX Titan X,
with CC values of 3.5 and 5.2, respectively. We note that the major
revision numbers are different. They correspond to two distinct
architectures, Kepler and Maxwell, with major revision numbers 3 and
5, respectively. These GPU devices are connected to two different
computers, having CPUs Dual-Core AMD Opteron 2222 working at a
frequency of 3\,GHz and Intel 7-4820K with a working frequency of
3.70\,GHz, respectively. We will refer to the CPUs of these systems as Kp and Mw,
respectively. In our first device we have 15 SMs, while in the second
one we have 24 SMs. We have decided to use both types of hardware
because, although the Maxwell microarchitecture is more recent and, a
priori, better in many technical aspects, it seems that double
precision arithmetics (required for practical applications of the CJM)
works better on Kepler than on
Maxwell. Table\,\mbox{\ref{tab:architectures}} summarizes the relevant
properties of the two architectures employed in this paper.
\begin{table}[ht!]
\centering
\begin{tabular}{|c c c c c|} 
\hline
CPU & Device                      &   Architecture &  CC  & SMs \\
\hline 
Opteron (Kp) & Tesla K40c               &   Kepler          &3.5  & 15  \\
Intel (Mx)     & GeForce GTX Titan X &  Maxwell        &5.2   & 24 \\
\hline
\end{tabular}
\caption{Properties of the two GPU architectures in which the CJM has
    been tested. We note that the memory on board (12 Gb) is the same
    for both GPUs.}
\label{tab:architectures}
\end{table}

\subsection{Code}
\label{subsec:code}

\begin{figure*}
\centering
\includegraphics[width=0.6\textwidth]{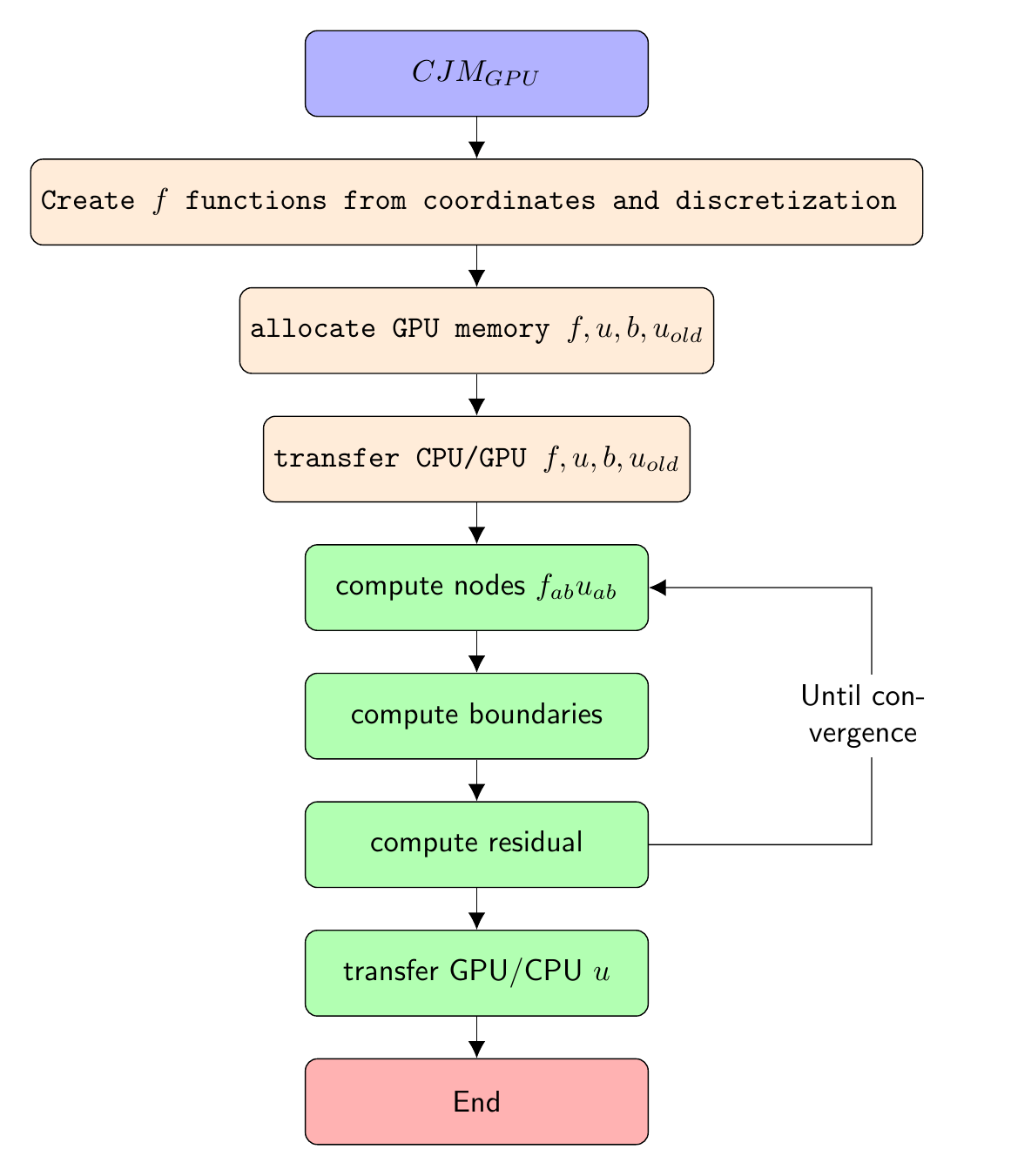}
\caption{Flow chart of the code. $b$ stands for the source term (if any) of the elliptic partial differential equation. Green background rectangles correspond to execution blocks in the GPU device. The rest of the execution blocks take place on the host computer.}
\end{figure*}

In order to execute the Jacobi method and the CJM over GPUs, we have
developed two CUDA kernels, i.e., sections of code that instead of
running on the CPU run on the GPU device. We need to transfer the data
that the kernel will use from the memory of the host to that of the
device. Once the results are computed, they are transferred back from
the device to the host. The 12GB of memory of our devices
is large enough to allow us to transfer the whole
data structure of the problem to the GPU memory in one transfer at the beginning
of the calculation. Once the problem has been solved
  to the desired accuracy (fully on the GPU device), we recover
back the solution in the CPU also with a single
data transfer. We have executed
the same code in the two considered GPU architectures.

\begin{table}[p]
\centering

\begin{tabular}{|c|c|c|}
\hline
$f_{\rm NW}(x_{1,i},x_{2,j},\Delta x_{1,i}, \Delta x_{2_j})$ & $f_{\rm N}(x_{1,i},x_{2,j},\Delta x_{1,i}, \Delta x_{2,j})$ & $f_{\rm NE}(x_{1,i},x_{2,j},\Delta x_{1,i}, \Delta x_{2,j})$ \\ \hline
$f_{\rm W}(x_{1,i},x_{2,j},\Delta x_{1,i}, \Delta x_{2,j})$ & $f_{\rm C}(x_{1,i},x_{2,j},\Delta x_{1,i}, \Delta x_{2,j})$ & $f_{\rm E}(x_{1,i},x_{2,j},\Delta x_{1,i}, \Delta x_{2,j})$ \\ \hline
$f_{\rm SW}(x_{1,i},x_{2,j},\Delta x_{1,i}, \Delta x_{2,j})$ & $f_{\rm S}(x_{1,i},x_{2,j},\Delta x_{1,i}, \Delta x_{2,j})$ & $f_{\rm SE}(x_{1,i},x_{2,j},\Delta x_{1,i}, \Delta x_{2,j})$ \\ \hline
\end{tabular}

\vspace{0.5cm}

\begin{tabular}{|c|c|c|c|c|}
\hline
$f_{\rm NW2}$  & $f_{\rm NWN}$ & $f_{\rm N2}$  & $f_{\rm NEN}$ &  $f_{\rm NE2}$ \\ \hline
$f_{\rm NWW}$ & $f_{\rm NW}$   & $f_{\rm N} $   & $f_{\rm NE}$    & $f_{\rm NEE}$ \\ \hline
$f_{\rm W2}$     & $f_{\rm W}$       & $f_{\rm C}$      & $f_{\rm E}$         & $f_{\rm E2}$    \\ \hline
$f_{\rm SWW}$ & $f_{\rm SW}$   & $f_{\rm S}$    & $f_{\rm SE}$    & $f_{\rm SEE}$  \\ \hline
$f_{\rm SW2}$  & $f_{\rm SWS}$ & $f_{\rm S2} $ & $f_{\rm SES}$ & $f_{\rm SE2}$   \\ \hline
\end{tabular}

\vspace{0.5cm}
\caption{Notation for the functions of neighbors at a distance one and two.}
\label{tab:ste1}
\end{table}

\begin{table}[p]
\centering

\begin{tabular}{|c|c|c|} 
\hline
\multicolumn{3}{|c|}{Cartesian coordinates}\\
\hline
$0$ & $\frac{1}{\Delta y^2}$ & $0$ \\ \hline
$\frac{1}{\Delta x^2}$ & $\frac{-2}{\Delta x^2} + \frac{-2}{\Delta y^2}$  & $\frac{1}{\Delta x^2}$ \\ \hline
$0$ & $\frac{1}{\Delta y^2}$  & $0$  \\ \hline
\end{tabular}

\vspace{0.5cm}

\begin{tabular}{|c|c|c|} \hline
\multicolumn{3}{|c|}{Polar coordinates}\\
\hline
$0$ & $\frac{1}{r_i^2 \Delta \theta^2}$ & $0$ \\ \hline
$\frac{1}{\Delta r_i^2} - \frac{1}{2 r_i \Delta r}$ & $\frac{-2}{\Delta r^2} + \frac{-2}{r_i^2 \Delta \theta^2}$  & $\frac{1}{\Delta r_i^2} + \frac{1}{2 r_i \Delta r}$ \\ \hline
$0$ & $\frac{1}{r_i^2 \Delta \theta^2}$  & $0$  \\ \hline
\end{tabular}

\vspace{0.5cm}

\begin{tabular}{|c|c|c|} \hline
\multicolumn{3}{|c|}{Bipolar coordinates}\\
\hline
$0$ & $\frac{\cosh \nu_j - \cos \mu_i}{a^2 \Delta \nu^2}$ & $0$ \\ \hline
$\frac{\cosh \nu_j - \cos \mu_i}{a^2 \Delta \mu^2}$ & $\frac{-2(\cosh \nu_j - \cos \mu_i)}{a^2 \Delta \mu^2}$ + $\frac{-2(\cosh \nu_j - \cos \mu_i)}{a^2 \Delta \nu^2}$  & $\frac{\cosh \nu_j - \cos \mu_i}{a^2 \Delta \mu^2}$ \\ \hline
$0$ & $\frac{\cosh \nu_j - \cos \mu_i}{a^2 \Delta \nu^2}$  & $0$  \\ \hline
\end{tabular}

\vspace{0.5cm}
\caption{Values of the variables of first table in Table
  \ref{tab:ste1} for some specific cases. In all cases we assume
  uniform meshes. The upper, middle and lower tables correspond to
    the standard 5-points discretization of the Laplacian operator in Cartesian, polar
    and bipolar coordinates, respectively.}
\label{tab:ste2}
\end{table}

The code is implemented as follows. The code is of stencil type, i.e.,
each node in the grid requires of knowing the values of the variable,
whose solution we seek, at a certain number of neighbors as
Fig.\,\mbox{\ref{fig::fig01}} schematically shows for the cases in
which the Laplacian operator is discretized in 9 or 17 points. In the
most general case we have implemented, both the central node,
identified with the integer indices $(i,j)$, where $1\le i \le N_x$
and $1\le j \le N_y$ and each of its (at most) 24 neighbors spanned by
the discretization of the Laplacian can have different numerical
factors weighting their contributions (see, e.g., the values of such
coefficients enclosed in the circles shown in
Fig.\,\mbox{\ref{fig::fig01}}, which correspond to the simplest case
of a two-dimensional problem in Cartesian coordinates). These
numerical factors may change as a function of the position of the
central node of the discretization in the grid. We refer to each of
these numerical factors with the terminology sketched in
Tab.\,\mbox{\ref{tab:ste1}}. The numerical factor corresponding to the
central node at a position $(i,j)$ is annotated by
$f_{\rm C}(x_{1,i},x_{2,j},\Delta x_{1,i}, \Delta x_{2,j})$, as it may
depend on the general coordinates $(x_{1,i},x_{2,j})$ and on the grid
spacing between consecutive nodes around the central node
$(\Delta x_{1,i},\Delta x_{2,j})$, since we allow for non-uniform
meshes. The rest of the numerical factors are named using the cardinal
points and numbers (see Tab.\,\ref{tab:ste1}). To improve
  the efficiency and to deal more easily with the boundary
  conditions, in the cases in which the Laplacian is discretized in
up to 9 points, we have created a different data structure for this
case (upper part of Tab.\,\mbox{\ref{tab:ste1}}), which we
differentiate from the most generic case represented by the lower part
of Tab.\,\mbox{\ref{tab:ste1}}.

As a result of this, the code is totally generic regarding
discretization and coordinates: the choice of coordinates determines
the discretization of the Laplacian operator, and this discretization
determines the mask of the functions (see Table \ref{tab:ste2}). In
addition, although the kernels are quite generic, the code allows
particular instantiations to the different stencils in a simple and
optimal way (without saving repeated values, without operating the
zeros, etc.).

Our implementation also takes advantage of memory hierarchy
on the GPU device, since the kernel uses both the Shared
Memory and the Registers on board of the CUDA
  device. Employing
  the CUDA Occupancy Calculator (free tool provided by Nvidia), we
  have tunned a number of parameters of the developed kernel in order
  to maximize the occupancy of the GPU devices. After some
  experimentation, we find that a 100\% occupancy results employing 
  256 threads per block and 2048 bytes of shared memory per block
  independent of the compute capability of the device. The number of
  registers per thread has been chosen to be 37 (32) for a device with
  CC\,$=3.5$ (CC\,$=5.2$).

\subsection{Test problem}
\label{subsec:tesPro}

We employ a very simple setting to
  calibrate the new implementation of our algorithm for GPUs. For that
  we solve a simple Poisson problem with a source term, whose solution
  is analytically known \mbox{\cite{Adsuara2017446}}:
\begin{equation}
\Delta u = -(x^2 + y^2) e^{x y}, \;\; (x,y)\in[0,1]^2,
\label{eq:Poisson2D}
\end{equation}
with appropriate Dirichlet boundary conditions. The unit
  square is discretized in $N_x\times N_y=1024\times 1024$ nodes,
  where the numerical solution is computed. The boundaries are easily
specified in this case, since there exists an analytic solution for
the problem at hand, which can be used to compute the boundaries at
the edges of the computational domain. The analytic solution reads
\begin{equation}
u(x,y)=-e^{x y}.
\label{eq:solPoisson2D}
\end{equation}
In next sections, we will take advantage of the knowledge of the analytic solution to compute the reduction of the error with resolution.

\subsection{Times and ratios}
\label{subsec:res}

We have solved the test problem of Sec.\,\ref{subsec:tesPro}, until
reaching a prescribed tolerance, using the classical 5-points
discretization of the Laplacian, and also the two discretizations of
9-points and 17-points introduced in
Sec.\,\ref{subsec:cjHigOrd}. We employ a resolution dependence
  tolerance, which decreases as $N^{-2}$.
For reference, we have solved the problem with the Chebyshev-Jacobi
method as well as with the classical Jacobi one. Finally, as we have
mentioned before in Sec.\,\ref{subsec:gpuDev}, we have used two
different microarchitectures for the implementations which work over
GPUs.

\begin{figure}[h!]
\centering
\includegraphics[width=0.49\textwidth]{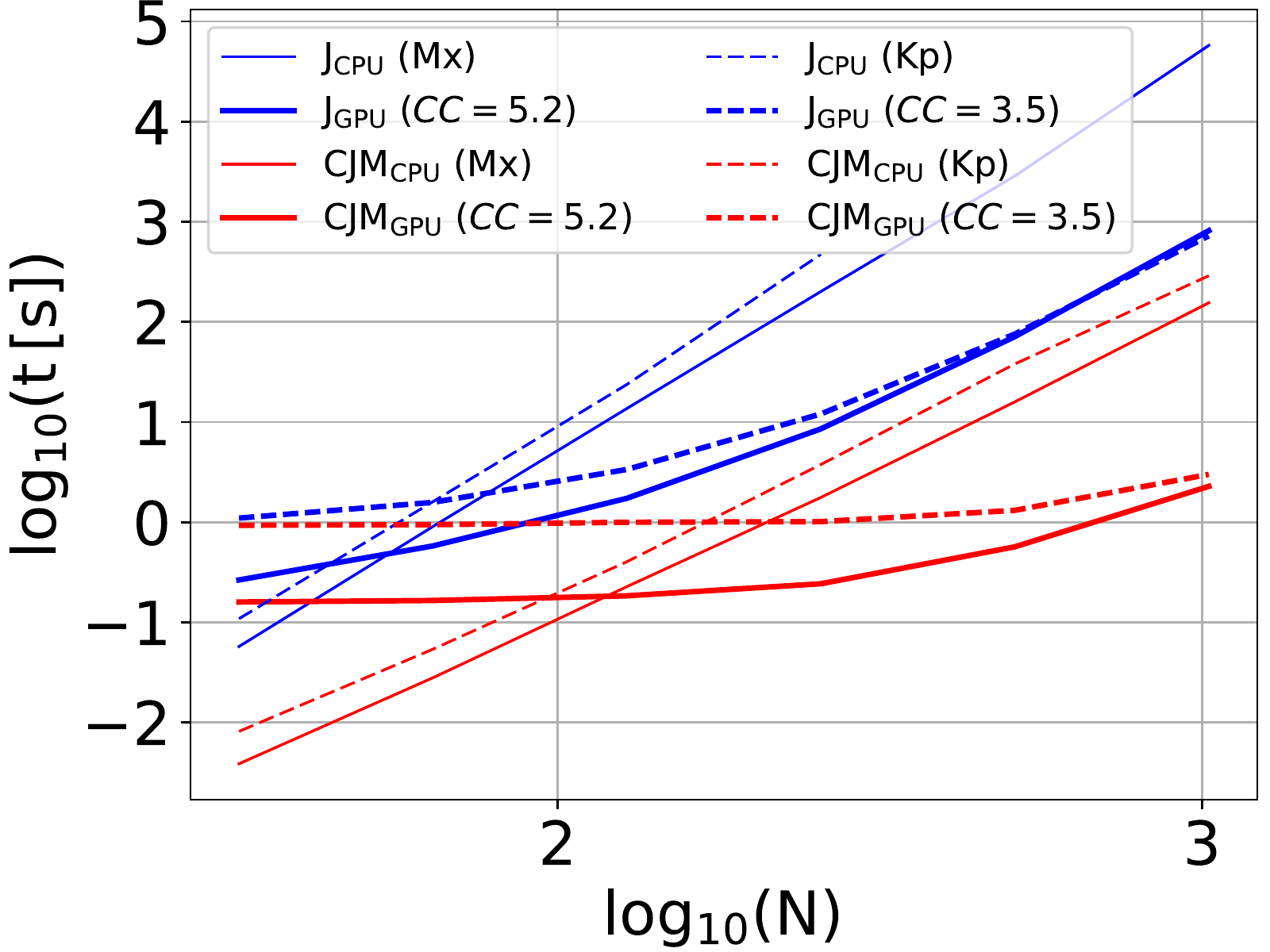}
\includegraphics[width=0.49\textwidth]{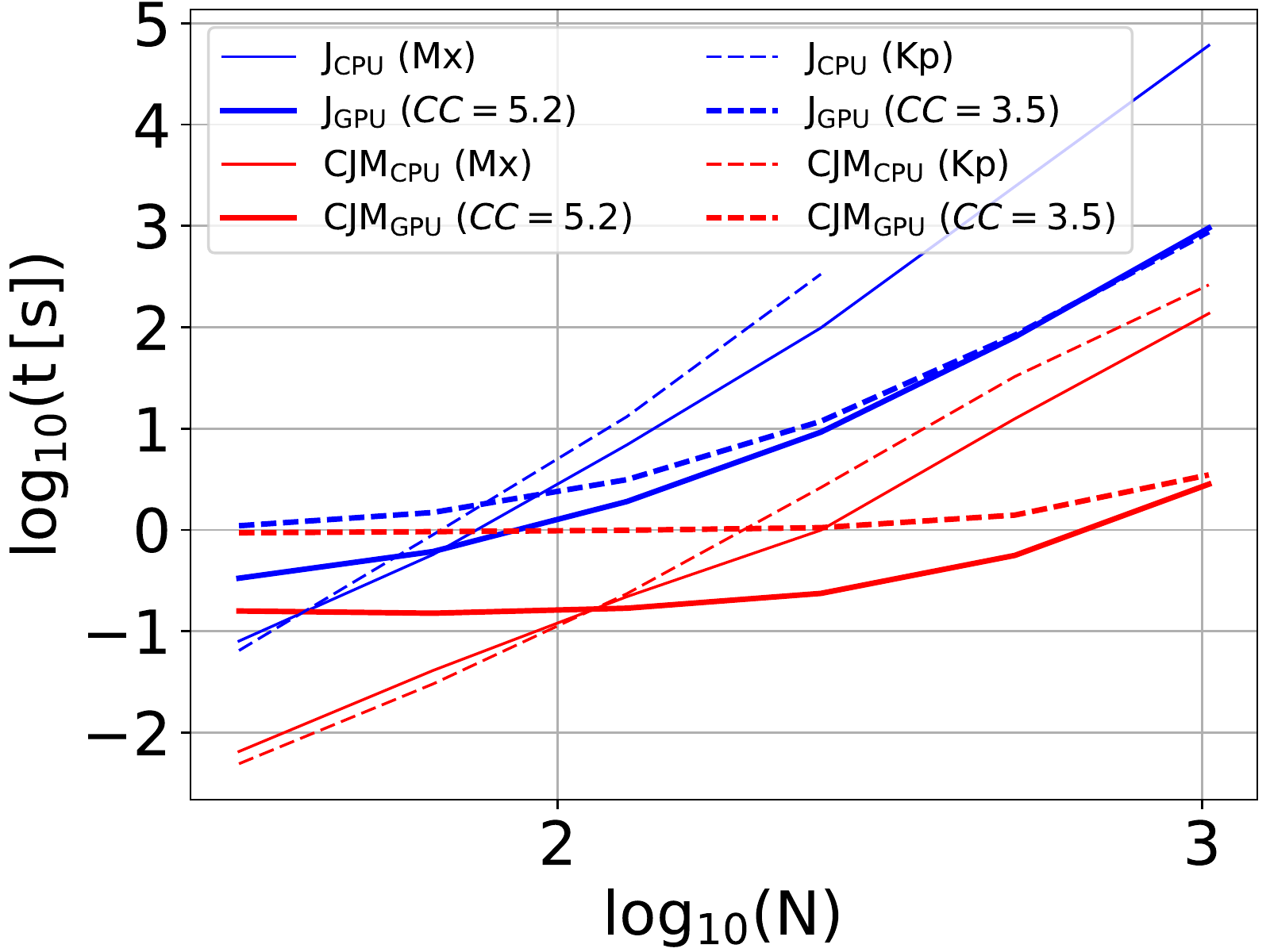}
\includegraphics[width=0.49\textwidth]{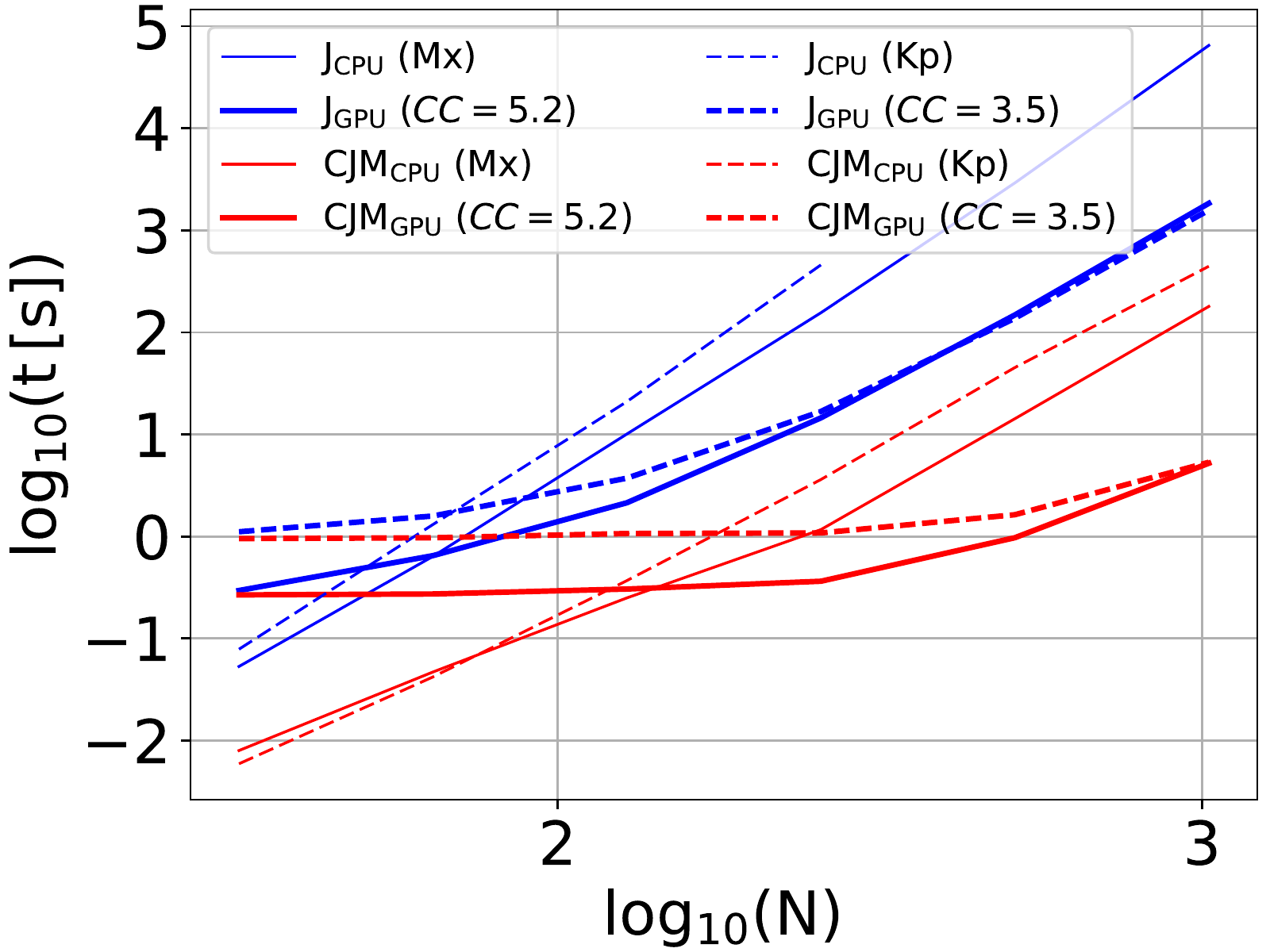} 
\includegraphics[width=0.49\textwidth]{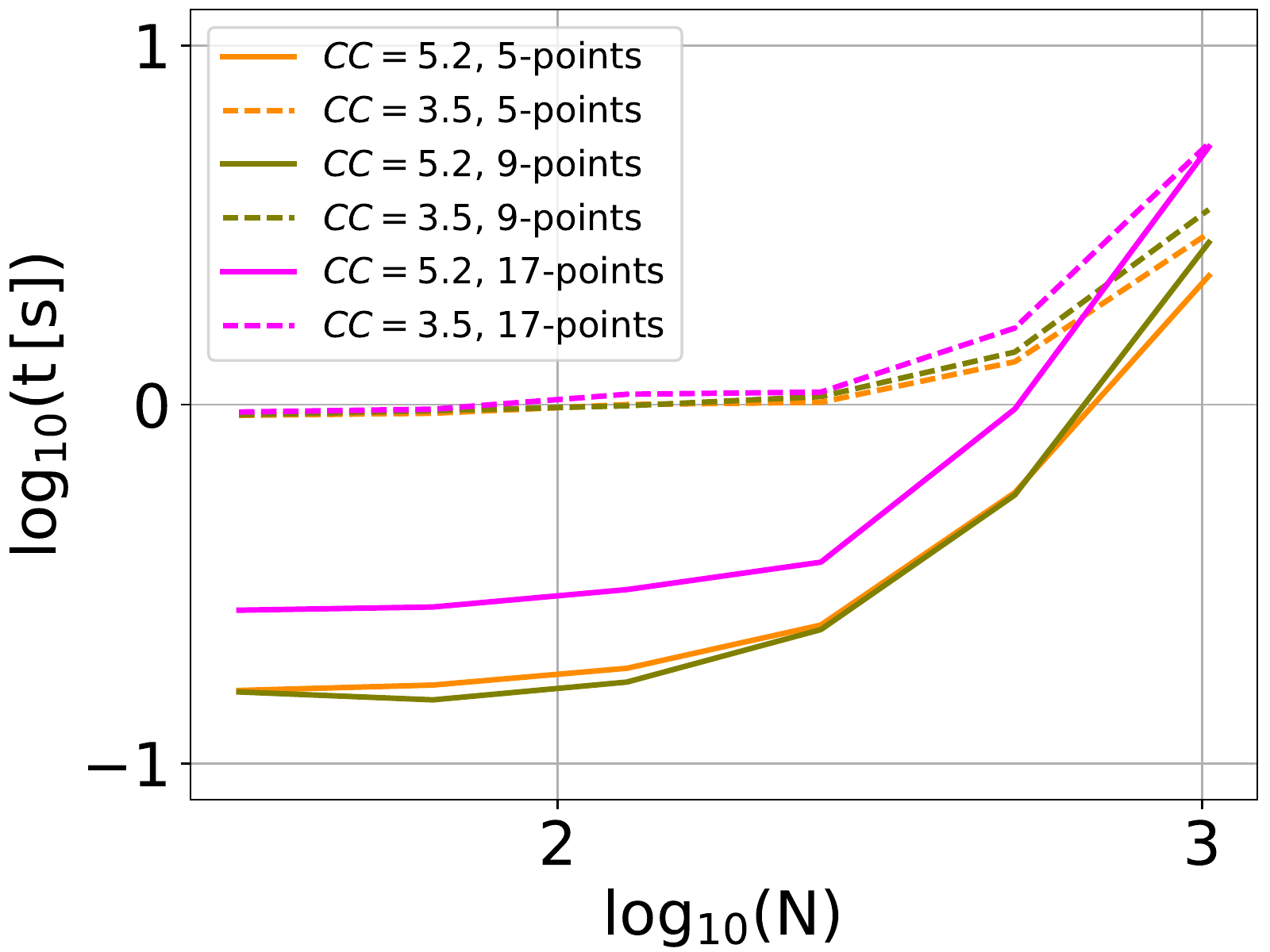} 
\caption{Both top and bottom left plots show the time necessary to
  reach a prescribed tolerance for different methods and
  implementations. We use a 5-points, 9-points and 17-points stencils
  in top left, top right and bottom left plots, respectively. Blue
  (red) lines refer to the Jacobi (Chebyshev-Jacobi) method. Thin
  (thick) lines correspond to serial code (GPU implementation). Solid
  (dashed) lines are associated to a Compute Capability equals to 5.2
  (3.4). Bottom right plot shows the times of the Chebyshev-Jacobi
  method over GPUs.}
\label{fig:fig03}
\end{figure}

Figure\,\mbox{\ref{fig:fig03}} shows the computational time (in
seconds) as a function of the number of points per dimension $N$. For
simplicity, we make our measurements in uniform computational meshes
satisfying $N_x=N_y=N$, so that the grid resolution in any of the two
spatial directions is $h=1/N$.  For low resolutions most of the time
in the GPU implementations is consumed by the transfer of data between
the host and the device. This time to transfer data is larger than the
time required by the method itself. In the different panels of
Fig.\,\mbox{\ref{fig:fig03}} this effect shows up as a plateau region
in the curves corresponding to GPU implementations of, specially, the
CJM. The plateau extends up to a certain turnover value $N_{\rm to}$
(depending on the method and on the architecture of the device), above
which we observe that the slopes of the lines stabilize in a very
similar way to the one of the sequential execution, but almost two
orders of magnitude below it. For instance, in the upper left panel,
of Fig.\,\ref{fig:fig03} this flat region extends up to
$N_{\rm to}\lesssim 300$. The latency of the data transfers is more
obvious in the Kepler architecture (corresponding to a CC of 3.5) than
in the case of Maxwell. The turnover in the former device happens at
$N_{\rm to}\sim 300$ for the CJM run over GPUs, while it is located at
$N_{\rm to}\sim 256$ for the latter device (see Fig.\,\ref{fig:fig03}
upper left panel). For the Jacobi method implemented on GPUs, the
aforementioned plateau does not exist. The reason for it is
  twofold. First, the amount of data transfers between the host and
  the CUDA device is slightly smaller in the latter method. This extra
  data required by the CJM are the $M$ weights needed to perform the $M$
  iterations of a complete computational cycle. Second, the Jacobi
  method is computationally more intensive than the CJM, since the the
  number of operations per grid node is comparable in both methods,
  but the number of iterations to reach the tolerance goal is much
  larger in the Jacobi method than in CJM. As a result, the Jacobi
  method is relatively more costly with respect to the data transfer.
We therefore conclude that a minimum mesh size is needed so that the
data transfer time between the CPU and the GPU, and viceversa, is
negligible in comparison to the computing time. Furthermore, we note
that below certain critical mesh size it is not even advantageous
using the CUDA devices, since the CPU implementation of the methods at
hand run faster. This is the case, e.g., of the Jacobi method when
using the 5-points discretization of the Laplacian displayed in
Fig.\,\ref{fig:fig03} (upper left panel): where the thin solid
(dashed) blue line, corresponding to the CPU executions, exhibit an
smaller computational time than the corresponding GPU runs, indicated
with thick solid (dashed) lines. Both lines (thin and thick
corresponding to the same method) cross at a ``critical'' value
$N_{\rm c}\lesssim 64$. This effect is exacerbated in the test
involving the CJM, where the CPU executions are faster than the
corresponding GPU ones for $N_{\rm c}\lesssim 128$. Noteworthy,
comparing the upper left panel of Fig.\,\mbox{\ref{fig:fig03}} with
the upper right and bottom left panels of the same figure, the
effectiveness of the GPU implementation of the CJM over the
corresponding CPU does not depend on the stencil of the discretization
of the Laplacian.

From the bottom right plot, where only the Chebyshev-Jacobi method on
GPUs is plotted, we observe that the qualitative behavior is
quite similar for all the different stencils. Comparing results obtained in the two CUDA devices with
  the same discretization of the Laplacian, it is evident that the
  architecture of Maxwell (CC\,$=5.2$) is faster than that of
  Kepler. However, the gap between both architectures reduces as the
  resolution increases. Actually, the difference in execution time
  between both devices displays a trend to reduce as the number of
  points in the discretization of the Laplacian grows, particularly
  beyond the value of $N_{\rm to}$ for each test. This is remarkably
  visible for the 17-points discretization of the Laplacian (pink
  lines of Fig.\,\mbox{\ref{fig:fig03}} bottom right panel), since at
  the highest resolution of our tests ($N=1024$), the computational
  time is the same in both devices.
This happens in spite of the fact that Maxwell
architecture is faster regarding both transfers and
executions. However, Kepler architecture is better suited
  for double precision computations, which we employ in practical
  applications of the elliptic solver at hand.

\begin{figure}[hp!]
\centering
\begin{tabular}{cc}
\hspace{-0.2cm}\includegraphics[width=0.49\textwidth]{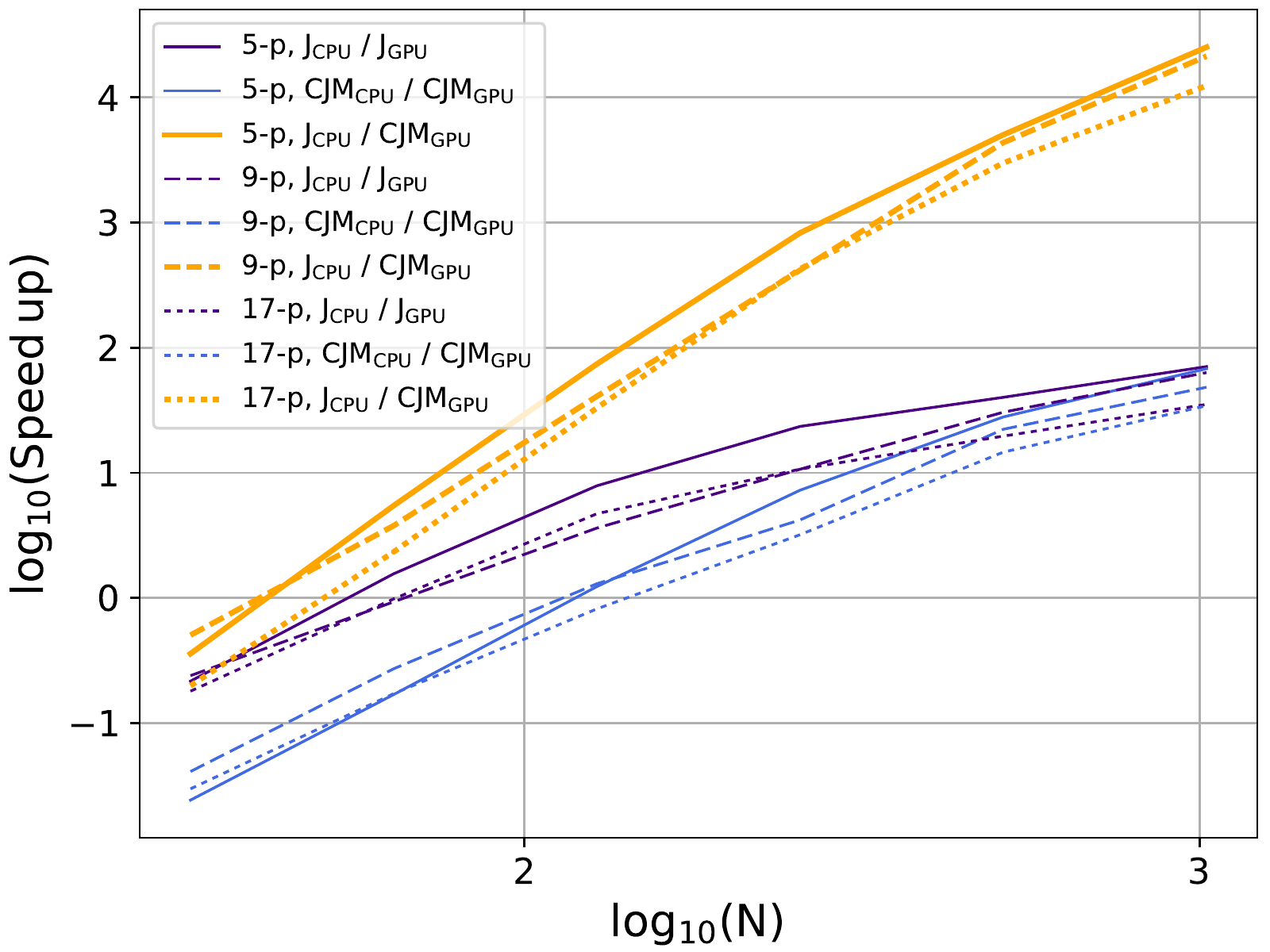} &\hspace{-0.55cm}
\includegraphics[width=0.49\textwidth]{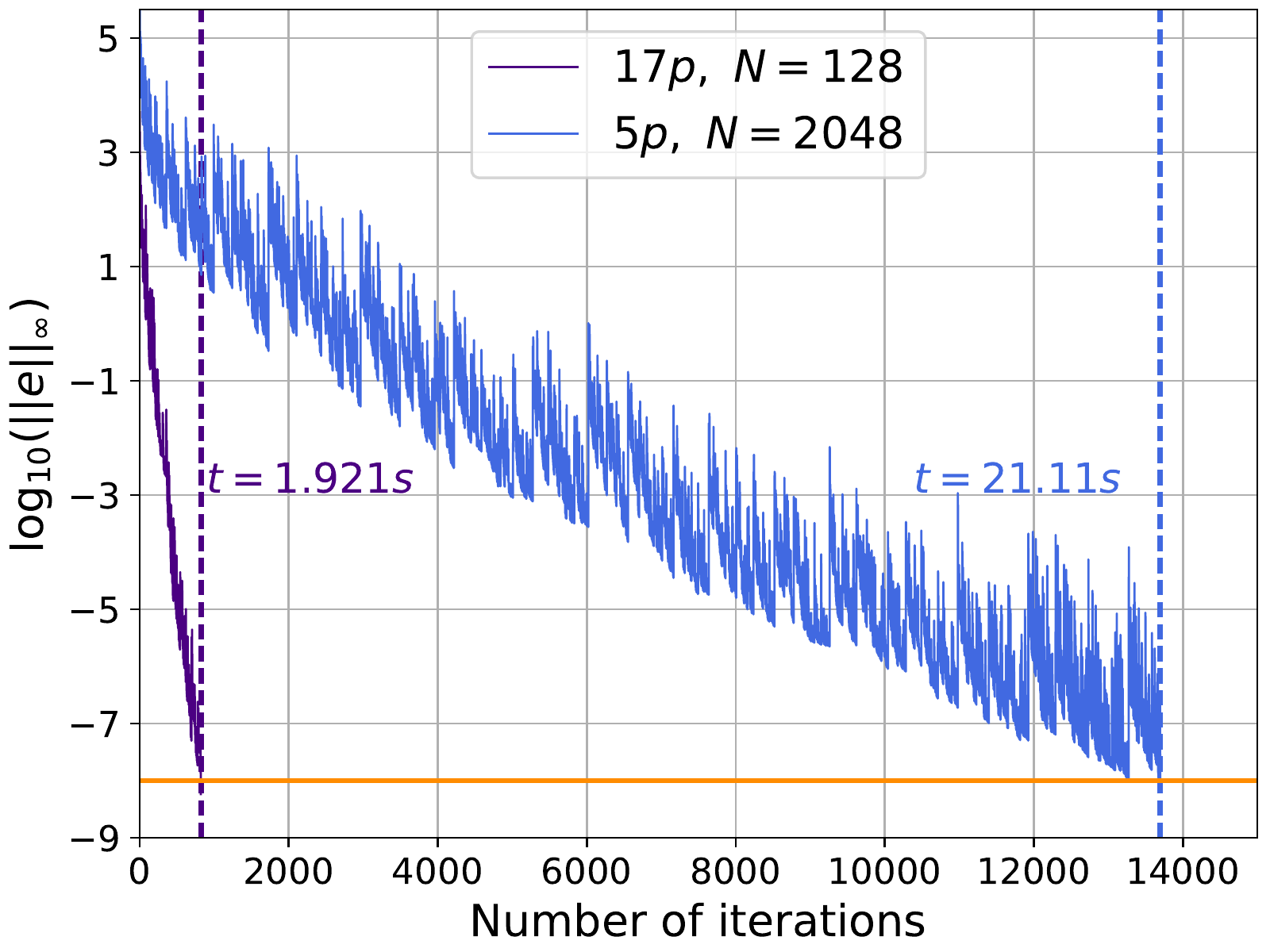} 
\end{tabular}

\caption{Left: Ratio of
    computational time in a CPU execution to the computational time in
    an analogous GPU execution as a function of the number of points
    per dimension. Solid lines correspond
  to the 5-points stencil, dashed lines to the 9-points stencil and
  dotted lines to the 17-points stencil. Blue and purples lines
  represent the ratios between the times of the sequential versus the
  parallel implementations of the Jacobi and Chebyshev-Jacobi
  methods. Orange lines represent the ratio between the time of the
  slowest method, the sequential implementation of the Jacobi method,
  and the fastest one, the Chebyshev-Jacobi method in parallel over
  GPUs. Right: Real error versus number of iterations for
  a 17-points stencil with a mesh of $N=128$
  points (purple line) and a classical 5-points stencil with a mesh of
  $N=2048$ points (blue line), until reaching a
  real error of $10^{-8}$ (orange horizontal line). The vertical
  dashed lines indicate the final number of iterations. The
  computational times of each run are also
  labeled in the plot.}
	\label{fig:fig04}
\end{figure}

In Figure \ref{fig:fig04} we plot the ratios
between the times of different methods and
implementations. These ratios provide an estimation of the
  acceleration factor by which the GPU implementation of either the
  CJM or the Jacobi method is faster (or slower) than any other of the
  methods. The figure also displays the
  dependence of the speed up of the GPU versions on the three types
of stencil. We observe that the larger number of
points the stencil has, the smaller is
  the speed up. The improvement over the CPU
  implementation is larger for the Jacobi method than for the CJM
  (note that the purple solid line lies above any of the blue lines in
  left panel of Fig.\,\mbox{\ref{fig:fig04}}). However, this larger
  relative speed up in the Jacobi method, with the smallest number of
  discretization points of the Laplacian, reduces with either
  increasing values of $N$, or when considering higher-order
  discretizations of the Laplacian. In the latter case, the speed up
  of the CJM over GPUs equals that of the corresponding Jacobi method
  employing a 17-points discretization of the Laplacian when
  $N=1024$. In fact, extrapolating the results to even higher resolutions
(not included in the plot) we foresee that the speed up
  factor of the CJM over its corresponding sequential CPU
  implementation will be better than that of the Jacobi
  method. From the ratio between the
CJM on GPUs, and
the slowest one, the classical Jacobi method in sequential, it can be
appreciated the difference of several orders of magnitude in the speed up
for high resolutions (see orange lines in left panel of
  Fig.\,\mbox{\ref{fig:fig04}}).

\begin{table}
\centering
\begin{tabular}{c|cccc}
    \hline
    $\texttt{5-points}$ & $\texttt{j}$ & $\texttt{j\_GPU}$ & $\texttt{cj}$ & 
$\texttt{cj\_GPU}$  \\ \hline
    $\texttt{j}$ & $\colorbox{white}{1}$  & $-$ & $-$ & $-$ \\ 
    $\texttt{j\_GPU}$ & $\colorbox{white}{71} $  & $\colorbox{white}{1}$ & $-$ & $-$ \\ 
    $\texttt{cj}$ & $\colorbox{white}{371}$ & $\colorbox{white}{5} $ & $\colorbox{white}{1}$ & 
$-$ \\
    $\texttt{cj\_GPU}$ & $\colorbox{white}{25235} $ & $\colorbox{white}{357}  
$ & $\colorbox{white}{68}$ & $\colorbox{white}{1}$ \\  
    \hline
    $\texttt{9-points}$ & $\texttt{j}$ & $\texttt{j\_GPU}$ & $\texttt{cj}$ & 
$\texttt{cj\_GPU}$  \\ \hline
    $\texttt{j}$ & $\colorbox{white}{1}$  & $-$ & $-$ & $-$ \\ 
    $\texttt{j\_GPU}$ & $\colorbox{white}{63}$  & $\colorbox{white}{1}$ & $-$ & $-$ \\ 
    $\texttt{cj}$ & $\colorbox{white}{441} $ & $\colorbox{white}{7} $ & $\colorbox{white}{1}$ & 
$-$ \\
    $\texttt{cj\_GPU}$ & $\colorbox{white}{21360} $ & $\colorbox{white}{337} $ 
& $\colorbox{white}{48}$ & $\colorbox{white}{1}$ \\ 
    \hline
    $\texttt{17-points}$ & $\texttt{j}$ & $\texttt{j\_GPU}$ & $\texttt{cj}$ & 
$\texttt{cj\_GPU}$  \\ \hline
    $\texttt{j}$ & $\colorbox{white}{1}$  & $-$ & $-$ & $-$ \\ 
    $\texttt{j\_GPU}$ & $\colorbox{white}{35} $  & $\colorbox{white}{1}$ & $-$ & $-$ \\ 
    $\texttt{cj}$ & $\colorbox{white}{361} $ & $\colorbox{white}{10} $ & $\colorbox{white}{1}$ & 
$-$ \\
    $\texttt{cj\_GPU}$ & $\colorbox{white}{12420} $ & $\colorbox{white}{352} $ 
& $\colorbox{white}{34} $ & $\colorbox{white}{1}$ \\
    \hline
  \end{tabular}
  \caption{Speed up factors of each of the methods listed in the
      columns with respect to the methods annotated in the rows in which
      the stencil of the Laplacian discretization is provided. All the
      values correspond to the test problem of Sec.\,\mbox{\ref{subsec:tesPro}}
      evaluated on a grid with $N=1024$ points per dimension.}
\label{tab:tab01}
\end{table}

In addition to Fig.\,\ref{fig:fig04}, we list the values
of the ratios for the particular case of a mesh with
$N=1024$ points per dimension in Table \ref{tab:tab01}. In this table we systematically
  include the speed up factors of each of the methods listed in the
  columns with respect to the methods annotated in the rows in which
  the stencil of the Laplacian discretization is provided.

Finally, as we have already shown in \cite{Adsuara2017446}, it proves
convenient to use a higher order discretization of the operator
although it involves larger
  stencils. In the latter case the computational time due to the fact
of having more points to get a solution of the
same quality is negligible compared to the reduction in time due to
the smaller number of iterations needed. To illustrate this point, we
use the CJM on GPUs (our
method of choice). We
  understand that having solutions of the same quality means
reaching the same ''real''
error level.  By real error we denote the difference, in
infinity norm, between the obtained numerical solution and the
analytical one. Since we have the analytical solution
\eqref{eq:solPoisson2D} of the the test problem \eqref{eq:Poisson2D},
we can use the real error instead of a tolerance as the stopping
criterion. In Fig.\,\ref{fig:fig04} we plot the real error versus
number of iterations using, on one hand, a 17-points stencil with a
mesh of $N=128$ points and, on the other hand, a
classical 5-points stencil with a mesh of
$N=2048$ points, until reaching a prescribed
value of the real error ($10^{-8}$); in addition to the smaller
resolution needed for the higher order method (17-points stencil), we
have a reduction of one order of magnitude both in the
number of iterations and in the computational
time.

\section{Conclusions}
\label{sec:con}

In previous papers \cite{
    Adsuara2016369,Adsuara2017446}, we have delineated the potencial
  degree of parallelism of the CJM method as a major advantage over
  other competing algorithms to solve linear systems of equations
  resulting from the discretization of elliptic systems of partial
  differential equations. Building upon the basic Jacobi method, the
  parallel implementation of the CJM is as simple as that of the
  former method. Here we have materialized our previous claims and
  presented a GPU based implementation of the CJM. We have tested the
  CUDA ported Jacobi and CJ algorithms in two different GPU
  architectures with compute capabilites of 3.5 and 5.2. Even though
  the Maxwell architecture (with compute capability 5.2) is more
  recent and potentially faster than the Kepler one (with CC\,$=3.5$),
  the differences in actual computing time reduce significantly as
  either the grid size increases or the number of points employed in
  the discretization of the Laplacian grows. We find that it is
possible to speed up by several orders of magnitude the classical Jacobi
method thanks to the use of the parallel implementation of the CJM on GPUs.

Moreover, we have illustrated the benefits from using the parallel
implementation of the CJM over
GPUs in combination with a high-order discretization of the Laplacian
operator in a test problem. We conclude that it is always
  advantageous employing high-order discretizations of the elliptic
  operator since they requiere less iterations and less computational
  time to reach the same real error goal. This conclusion is
  independent of the parallel implementation of any of the methods we
  have tested in this paper. However, the combination of high-order
  discretization of the elliptic operators and the CJM implemented on
  GPUs results in an extremely powerful method for practical
  applications.

\section*{Acknowledgments}
We thank Marc Jord\`{a} from the Barcelona Supercomputing Center for his contribution in the implementations using GPUs. We acknowledge the support of the Consolider SyeC net TIN2014-52608-REDC. We also acknowledge the support from the European Research Council (Starting Independent Researcher Grant CAMAP-259276) and from the Spanish Ministerio de Econom\'{\i}a y Competitividad through the grants SAF2013-49284-EXP, AYA2015-66899-C2-1-P, as well as the support of the Valencia Government through the grant PROMETEO-II-2014-069.


\bibliographystyle{elsarticle-num} 
\bibliography{cjGPUs}

\begin{thebibliography}{1}
\expandafter\ifx\csname url\endcsname\relax
  \def\url#1{\texttt{#1}}\fi
\expandafter\ifx\csname urlprefix\endcsname\relax\def\urlprefix{URL }\fi
\expandafter\ifx\csname href\endcsname\relax
  \def\href#1#2{#2} \def\path#1{#1}\fi

\bibitem{YangMittal2014}
X.~I. Yang, R.~Mittal,
  \href{http://www.sciencedirect.com/science/article/pii/S0021999114004173}{Acceleration
  of the jacobi iterative method by factors exceeding 100 using scheduled
  relaxation}, Journal of Computational Physics 274 (2014) 695 -- 708.
\newblock \href {http://dx.doi.org/http://dx.doi.org/10.1016/j.jcp.2014.06.010}
  {\path{doi:http://dx.doi.org/10.1016/j.jcp.2014.06.010}}.
\newline\urlprefix\url{http://www.sciencedirect.com/science/article/pii/S0021999114004173}

\bibitem{Adsuetal15}
J.~E. {Adsuara}, I.~{Cordero-Carri\'on}, P.~{Cerd\'a-Dur\'an}, M.~A. {Aloy},
  {Improvements in the Scheduled Relaxation Jacobi method}, in: C.~{D{\'{\i}az
  Moreno}, J.~M. and {D{\'{\i}}az Moreno}, J.~C. and {Garc{\'{\i}}a
  V{\'a}zquez}, C. and {Medina Moreno}, J. and {Orteg{\'o}n Gallego}, F.
  {P{\'e}rez Mart\'{\i}}nez}, M.~V. {Redondo Neble}, J.~R. {Rodr{\'{\i}}guez
  Galv{\'a}n} (Eds.), Proceedings of the XXIV Congress on Differential
  Equations and Applications / XIV Congress on Applied Mathematics C{\'a}diz,
  June 8-12, 2015., Servicio de Publicaciones de la Universidad de C{\'a}diz,
  2015.

\bibitem{Adsuara2016369}
J.~Adsuara, I.~Cordero-Carri\'on, P.~Cerd\'a-Dur\'an, M.~Aloy,
  \href{http://www.sciencedirect.com/science/article/pii/S002199911630198X}{Scheduled
  relaxation jacobi method: Improvements and applications}, Journal of
  Computational Physics 321 (2016) 369 -- 413.
\newblock \href {http://dx.doi.org/http://dx.doi.org/10.1016/j.jcp.2016.05.053}
  {\path{doi:http://dx.doi.org/10.1016/j.jcp.2016.05.053}}.
\newline\urlprefix\url{http://www.sciencedirect.com/science/article/pii/S002199911630198X}

\bibitem{Adsuara2017446}
J.~Adsuara, I.~Cordero-Carri\'on, P.~Cerd\'a-Dur\'an, V.~Mewes, M.~Aloy,
  \href{http://www.sciencedirect.com/science/article/pii/S0021999116306738}{On
  the equivalence between the scheduled relaxation jacobi method and
  richardson's non-stationary method}, Journal of Computational Physics 332
  (2017) 446 -- 460.
\newblock \href {http://dx.doi.org/http://dx.doi.org/10.1016/j.jcp.2016.12.020}
  {\path{doi:http://dx.doi.org/10.1016/j.jcp.2016.12.020}}.
\newline\urlprefix\url{http://www.sciencedirect.com/science/article/pii/S0021999116306738}

\bibitem{Adams:1988}
L.~M. {Adams}, R.-J. {LeVeque}, D.~{Young},
  \href{http://dx.doi.org/10.1137/0725066}{Analysis of the sor iteration for
  the 9-point laplacian}, SIAM Journal on Numerical Analysis 25~(5) (1988)
  1156--1180.
\newblock \href {http://dx.doi.org/10.1137/0725066}
  {\path{doi:10.1137/0725066}}.
\newline\urlprefix\url{http://dx.doi.org/10.1137/0725066}

\bibitem{CUDA}
\href{http://www.nvidia.com/object/cuda\_home\_new.html}{Cuda}.
\newline\urlprefix\url{http://www.nvidia.com/object/cuda\_home\_new.html}

\end{thebibliography}

\end{document}